\documentclass{amsart}

\usepackage{amsmath,amssymb,amsthm}
\usepackage{graphicx,tikz}

\newtheorem*{thm}{Theorem}

\theoremstyle{definition}

\theoremstyle{remark}

\DeclareMathOperator{\med}{med}

\begin{document}

\title[]{A Sharp estimate for probability distributions}
\keywords{Distribution-free estimate, comparison inequalities, geometry of probability distributions}
\subjclass[2010]{60E05, 60E15} 

\author[]{Stefan Steinerberger}
\address{Department of Mathematics, Yale University, New Haven, CT 06511, USA}
\email{stefan.steinerberger@yale.edu}

\begin{abstract}  We consider absolutely continuous probability distributions $f(x)dx$ on $\mathbb{R}_{\geq 0}$.
A result of Feldheim \& Feldheim shows, among other things, that if the distribution is not compactly supported, then there exist $z > 0$ such that most
events in $\left\{X + Y \geq 2z\right\}$ are comprised of a 'small' term satisfying $\min(X,Y) \leq z$ and a 'large' term satisfying $\max(X,Y) \geq z$ (as opposed
to two 'large' terms that are both larger than $z$)
$$ \limsup_{z \rightarrow \infty}~ \mathbb{P}\left( \min(X,Y) \leq z | X+Y \geq 2z\right) = 1.$$
The result fails if the distribution is compactly supported. We prove
$$\sup_{z > 0 } ~\mathbb{P}\left( \min(X,Y) \leq z | X+Y \geq 2z\right)  \geq  \frac{1}{24 + 8\log_2{( \med(X) \|f\|_{L^{\infty}})}},$$
where $\med X$ denotes the median. Interestingly, the logarithm is necessary and the result is sharp up to constants; we also discuss some open problems.
\end{abstract}

\maketitle

\section{Introduction and main result}
\subsection{Introduction.} The purpose of this short note is to describe a sharp estimate for absolutely continuous probability distributions. We start by discussing related results
in a similar spirit. The 123 Theorem of Alon \& Yuster \cite{alon}, originally suggested by G. Margulis, is named after the three constants: if $X, Y$
are two i.i.d. real random variables, then there is the beautiful inequality
$$ \mathbb{P}\left( |X-Y| \leq 2\right) \leq 3~ \mathbb{P}\left(|X-Y| \leq 1\right).$$
The paper of Alon \& Yuster \cite{alon} proves both the 123 inequality and establishes the more general result (with sharp constant)
$$ \mathbb{P}\left( |X-Y| \leq b\right) \leq (2 \left\lceil b/a \right\rceil + 1)~ \mathbb{P}\left(|X-Y| \leq a\right).$$
Another inequality in this spirit was given by Siegmung-Schultze \& von Weizs\"acker \cite{schultze} under the same assumptions: for any $c>0$
$$ \mathbb{P}\left( |X+Y| \leq c \right) \leq 2 ~\mathbb{P}\left(|X-Y| \leq c\right).$$
Even though these results are remarkable and of great intrinsic interest ('the unavoidable geometry of probability distributions'), there seems to be relatively
little work on problems of these type. Our result is inspired by a beautiful Theorem of Feldheim \& Feldheim \cite{feldheim} (answering a question of Alon): for any two independent nonnegative
random variables $X,Y$ with non-compactly supported distributions 
$$ \limsup_{z \rightarrow \infty}~ \mathbb{P}\left( \min(X,Y) \leq z | X+Y \geq 2z\right) = 1.$$
Put differently, the result says that there are arbitrarily large values $z > 0$ such that a typical
event in $\left\{ X + Y \geq 2z\right\}$ is create by a 'small' value $\min(X,Y) \leq z$ and a 'large' value $\max(X,Y) \geq z$ as opposed
to two 'large' terms both larger than $z$.
We emphasize that in this result $X,Y$ may come from different probability distributions; Feldheim \& Feldheim \cite{feldheim}
obtain even more general results for random variables coming from the same distribution.

\subsection{The Result.} We will allow for the distribution to be compactly supported but will require that  $X$ and $Y$ are coming from the same probability distribution
which is assumed to be absolutely continuous and given by a measure $f(x) dx$. We show that there is a value $z > 0$ such that a certain proportion of events in $\left\{X+Y \geq 2z\right\}$ is comprised of events having a small term $X \leq z$ and a large term $Y \geq z$.

\begin{thm} If $X,Y$ are i.i.d. random variables drawn from an absolutely continuous probability distribution with density $f(x)dx$ on $\mathbb{R}_{\geq 0}$, then
$$\sup_{z>0} ~\mathbb{P}\left( X \leq z ~ \mbox{and} ~ X+Y \geq 2z\right) \geq  \frac{1}{24 + 8\log_2{( \med(X) \|f\|_{L^{\infty}})}},$$
where $\med X$ denotes the median of the probability distribution. This estimate is sharp up to constants
and the supremum can be restricted to $0 \leq z \leq \med(X)$. 
\end{thm}
It seems interesting that such an elementary question in probability theory has such a nontrivial logarithmic dependence on the distribution.
Absolutely continuous probability distribution have a symmetry under dilation: the left-hand side is invariant under this symmetry and so is
$\med(X) \|f\|_{L^{\infty}}$. Note that
$$ \|f\|_{L^{\infty}} \geq \frac{1}{\med(X)}\int_{0}^{\med(X)}{f(x)dx} = \frac{1}{2\med(X)} \quad \mbox{and thus} \quad  \med(X) \|f\|_{L^{\infty}} \geq \frac{1}{2}$$
and therefore the right-hand side is always defined. Trivially,
$$ \mathbb{P}\left( X \leq z ~ \mbox{and} ~ X+Y \geq 2z\right) \leq \mathbb{P}\left( \min(X,Y) \leq z | X+Y \geq 2z\right)$$
which implies the result stated in the abstract. Moreover, since it is possible to restrict the range of $z$ in the
main result, we can also obtain a sharp result for the other quantity:  for $z \leq \med(X)$, 
$$ \mathbb{P}(X+Y \geq 2z) \geq \mathbb{P}(X+Y \geq 2\med(X)) \geq \mathbb{P}(X \geq \med X)^2 = \frac{1}{4}$$
and therefore the inequality
$$\max_{ 0 \leq z \leq \med X } ~\mathbb{P}\left( \min(X,Y) \leq z | X+Y \geq 2z\right)  \geq  \frac{1}{24 + 8\log_2{( \med(X) \|f\|_{L^{\infty}})}}$$
is also sharp up to (possibly different)  constants. There are many natural questions.
Is it possible to obtain similar quantitative estimates if $X$ and $Y$ are drawn from different probability distributions?

\subsection{An open problem.}  We can rewrite the
Feldheim \& Feldheim result as
$$ \limsup_{z \rightarrow \infty} \frac{ \mathbb{P}(X+Y \geq 2z  \quad \mbox{and} \quad \min(X,Y) \leq z) }{\mathbb{P}\left(X +Y \geq 2z  \quad \mbox{and} \quad \min(X,Y) \geq z\right)} = \infty.$$
An improved result would be one where one adds an additional weight to the denominator and still obtains unboundedness.
A computation shows that if the random variables $X,Y$ are distributed with the same exponential distribution then
$$ \forall~z \geq 0 \qquad \frac{ \mathbb{P}(X+Y \geq 2z  \quad \mbox{and} \quad \min(X,Y) \leq z) }{\mathbb{P}\left(X +Y \geq 2z \quad \mbox{and} \quad \min(X,Y) \geq z\right)} \frac{\mathbb{E} X}{z} = 2.$$
Are exponential distributions the sharp limiting case? There are two questions attached to this: one for distributions with median fixed and one for distributions with expectation fixed (the first obviously being limited to probability distributions whose mean exists).
 More precisely, is it true that for any $X,Y$ drawn from an absolutely continuous probability distribution on $\mathbb{R}_{\geq 0}$ that is not compactly supported
\begin{align*}
\exists~z>0 \qquad\frac{ \mathbb{P}(X+Y \geq 2z \quad \mbox{and} \quad \min(X,Y) \leq z) }{\mathbb{P}\left(X +Y \geq 2z \quad \mbox{and} \quad \min(X,Y) \geq z\right)}&\geq \left\{  \frac{(2 \log{2}) z}{\med(X) },   \frac{2z}{\mathbb{E}X } \right\} ?
\end{align*}
If exponential distributions are not extremal, does the result nonetheless hold with some other constants?

\section{proof}
\begin{proof} We wish to bound
$$ \delta := \max_{0 \leq z \leq \med(X)} \mathbb{P}\left( X \leq z ~ \mbox{and} ~ X+Y \geq 2z\right)$$
from below. This quantity is invariant under the dilation symmetry
$$ f(x)dx \rightarrow \lambda f(\lambda x) dx \qquad \mbox{for all}~\lambda > 0$$
mapping probability distribution on $\mathbb{R}_{\geq 0}$ to probability distributions on $\mathbb{R}_{\geq 0}$
which allows us to assume w.l.o.g. that $\med X = 1$. 
We can rewrite the probability as
\begin{align*}
\mathbb{P}\left( X \leq z ~ \mbox{and} ~ X+Y \geq 2z\right) &= \int_{0}^{z}{f(x) \int_{2z - x}^{\infty}{f(y) dy} dx} \\
&= \int_{0}^{z}{f(x)\left(1 - F(2z-x)\right) dx}, \end{align*}
 where $F$ denotes the cumulative distribution function.
Let us now assume that $1/2 \leq z \leq 1$ which implies $0 \leq 2z-1 \leq z$. We know that $F$ is monotonically increasing and $F(1) = 1/2$. Thus
\begin{align*}
 \delta &\geq  \int_{0}^{z}{f(x)\left(1 - F(2z-x)\right) dx} \\
 &\geq \int_{2z-1}^{z}{f(x)\left(1 - F(2z-x)\right) dx} \geq \frac{1}{2}\int_{2z-1}^{z}{f(x)dx}.
 \end{align*}
Setting $z = 1 - 1/2^k$ for $k\geq 1$ yields
$$ \int_{1-\frac{1}{2^{k-1}}}^{1-\frac{1}{2^{k}}}{f(x)dx} \leq 2\delta.$$
We can now estimate
\begin{align*}
 \frac{1}{2} = \int_{0}^{1}{f(x)dx} &\leq  \left(\sum_{k=1}^{\ell}{ \int_{1-\frac{1}{2^{k-1}}}^{1-\frac{1}{2^{k}}}{f(x)dx}}\right) + \int_{1-2^{-\ell}}^{1}{f(x)dx} \\
&\leq \sum_{k=1}^{\ell}{ \int_{1-2^{-k+1}}^{1-2^{-k}}{f(x)dx}} + 2^{-\ell} \|f\|_{L^{\infty}} \leq 2 \delta \ell +  2^{-\ell} \|f\|_{L^{\infty}}.
\end{align*}
Setting $\ell = \left\lfloor 3 + \log_2{ \|f\|_{L^{\infty}}} \right\rfloor$ implies $\delta \geq 1/(8\ell)$ and thus
$$ \delta \geq \frac{1}{24 + 8 \log_2{\|f\|_{L^{\infty}}}}.$$
Using the dilation symmetry yields the result in the case $\med(X) \neq 1$. It remains to show that the result is sharp up to constants. We consider the distribution
$$   \left(\sum_{k=1}^{n}{ \frac{50 \cdot 2^k}{n}  \chi_{[1-\frac{1.01}{2^k}, 1 - \frac{0.99}{2^k}]}    }\right)dx \qquad \mbox{satisfying} \qquad \|f\|_{L^{\infty}} = 50\frac{2^n}{n}$$
and thus $\log_2{\|f\|_{L^{\infty}}} = n + o(n)$. Moreover, $\med(X) \sim 1.$
Let us assume that $$ 1 - \frac{1}{2^{k}} < z \leq 1 - \frac{1}{2^{k+1}} \qquad \mbox{for some}~k \in \mathbb{N}.$$
Clearly, $\mathbb{P}(X \leq 1/2) = 1/(2n)$ is small and it suffices to consider $k \geq 1$. 
The event 
$$X+Y \geq 2z \geq 2 - \frac{1}{2^{k-1}} \quad \mbox{requires} \quad  X \geq 1 - \frac{1}{2^{k-1}}$$
because $Y < 1$.
This implies
\begin{align*}
 \mathbb{P}\left( X \leq z ~ \mbox{and} ~ X+Y \geq 2z\right)  \leq  \mathbb{P}\left( 1 - \frac{1}{2^{k-1}}  \leq X \leq  1 - \frac{1}{2^{k+1}}\right) 
 \end{align*}
 and an explicit computation shows that
\begin{align*}
\mathbb{P}\left( 1 - \frac{1}{2^{k-1}}  \leq X \leq  1 - \frac{1}{2^{k+1}}\right) = \frac{2}{n} = \frac{2 + o(1)}{\log_2{\|f\|_{L^{\infty}}}}.
 \end{align*}
The construction shows that any bound of the form
$$\sup_{z>0} ~\mathbb{P}\left( X \leq z ~ \mbox{and} ~ X+Y \geq 2z\right) \geq  \frac{1}{A + B\log_2{( \med(X) \|f\|_{L^{\infty}})}} ~  \mbox{has}~B \geq 1/2.$$
\end{proof}
\textbf{Acknowledgement.} This paper was written during a stay at the Wispl Institute; the author gratefully acknowledges its hospitality.


\begin{thebibliography}{10}
\bibitem{alon} N. Alon and R. Yuster, The 123 theorem and its extensions, J. Combin. Theory Ser. A 72, 322--331 (1995).

\bibitem{dong}  Z. Dong, J. Li, and W. Li, A note on distribution-free symmetrization inequalities, J. Theor.
Prob. 28(3) (2015), 958--967.

\bibitem{feldheim} N. Feldheim and O. Feldheim, Mean and minimum of independent random variables, arXiv: 1609.03004

\bibitem{schultze} R. S. Schultze and H. V. Weizs\"acker, Level-crossing probabilities I: One-dimensional random walks and
symmetrization, Adv. Math. 208, 672--679 (2007)
\end{thebibliography}
\end{document}